# A STABLE MARRIAGE OF POISSON AND LEBESGUE

By Christopher Hoffman,[1] Alexander E. Holroyd[2]
and Yuval Peres[3]

*University of Washington, University of British Columbia and UC Berkeley*

Let $\Xi$ be a discrete set in $\mathbb{R}^d$. Call the elements of $\Xi$ *centers*. The well-known Voronoi tessellation partitions $\mathbb{R}^d$ into polyhedral regions (of varying sizes) by allocating each site of $\mathbb{R}^d$ to the closest center. Here we study "fair" allocations of $\mathbb{R}^d$ to $\Xi$ in which the regions allocated to different centers have equal volumes.

We prove that if $\Xi$ is obtained from a translation-invariant point process, then there is a unique fair allocation which is *stable* in the sense of the Gale–Shapley marriage problem. (I.e., sites and centers both prefer to be allocated as close as possible, and an allocation is said to be *unstable* if some site and center both prefer each other over their current allocations.)

We show that the region allocated to each center $\xi$ is a union of finitely many bounded connected sets. However, in the case of a Poisson process, an infinite volume of sites are allocated to centers further away than $\xi$. We prove power law lower bounds on the allocation distance of a typical site. It is an open problem to prove *any* upper bound in $d > 1$.

**1. Introduction.** Let $d \geq 1$. We call the elements of $\mathbb{R}^d$ *sites*. We write $|\cdot|$ for the Euclidean norm and $\mathcal{L}$ for Lebesgue measure or *volume* on $\mathbb{R}^d$. Let $\Xi \subset \mathbb{R}^d$ be a discrete set. We call the elements of $\Xi$ *centers*. Let $\alpha \in [0, \infty]$ be a parameter, called the *appetite*. An *allocation* (of $\mathbb{R}^d$ to $\Xi$ with appetite $\alpha$) is a measurable function

$$\psi : \mathbb{R}^d \to \Xi \cup \{\infty, \Delta\}$$

Received June 2005; revised October 2005.
[1]Supported in part by MSRI, NSF Grant DMS-00-99814 and the University of Washington Research Royalty Fund.
[2]Supported in part by an NSERC (Canada) research grant and by MSRI, CPAM and PIMS.
[3]Supported in part by NSF Grants DMS-01-04073 and DMS-02-44479 and by MSRI and CPAM.
*AMS 2000 subject classification.* 60D05.
*Key words and phrases.* Stable marriage, point process, phase transition.

such that $\mathcal{L}[\psi^{-1}(\Delta)] = 0$, and

$$\mathcal{L}[\psi^{-1}(\xi)] \leq \alpha$$

for all $\xi \in \Xi$. We call $\psi^{-1}(\xi)$ the *territory* of the center $\xi$. We say that $\xi$ is *sated* if $\mathcal{L}[\psi^{-1}(\xi)] = \alpha$ and *unsated* otherwise. We say that a site $x$ is *claimed* if $\psi(x) \in \Xi$ and *unclaimed* if $\psi(x) = \infty$. Note that $\infty$ and $\Delta$ will have very different meanings: $\psi(x) = \infty$ means that $x$ is unable to find any center willing to accept it, whereas $\psi(x) = \Delta$ means that $\psi(x)$ is "undefined"—we allow an $\mathcal{L}$-null set of such sites purely as a technical convenience.

The following definition is an adaptation of that introduced by Gale and Shapley [6]. The idea is that sites and centers both prefer to be allocated as close as possible, and an allocation is said to be unstable if some site and center both prefer each other over their current allocations.

*Definition of stability.* Let $\xi$ be a center and let $x$ be a site with $\psi(x) \notin \{\xi, \Delta\}$. We say that $x$ *desires* $\xi$ if

$$|x - \xi| < |x - \psi(x)| \qquad \text{or } x \text{ is unclaimed.}$$

We say that $\xi$ *covets* $x$ if

$$|x - \xi| < |x' - \xi| \qquad \text{for some } x' \in \psi^{-1}(\xi), \text{ or } \xi \text{ is unsated.}$$

We say that a site-center pair $(x, \xi)$ is *unstable* for the allocation $\psi$ if $x$ desires $\xi$ and $\xi$ covets $x$. An allocation is *stable* if there are no unstable pairs. Note that no stable allocation may have both unclaimed sites and unsated centers.

Figures 1 and 2 illustrate stable allocations for several values of $\alpha$ on the 2-torus.

We consider three main questions about stable allocations. First, given a countable set in the plane, does there exist a stable allocation? Second, if there is a stable allocation, is it unique? Finally, what do the territories look like for a typical set of centers, qualitatively and quantitatively?

THEOREM 1 (Existence). *For any discrete set of centers $\Xi \subseteq \mathbb{R}^d$, and any appetite $\alpha \in [0, \infty]$, there exists a stable allocation.*

Note that some condition on $\Xi$ is needed to guarantee existence of a stable allocation. For example, if $d = 1$ and $\Xi = \{1, 1/2, 1/3, \ldots\}$, then it is easy to see that some center must form an unstable pair with the site $-1$.

We can establish uniqueness of the stable allocation if $\Xi$ is a sufficiently well-behaved set in two dimensions, or if $\Xi$ arises from a point process.

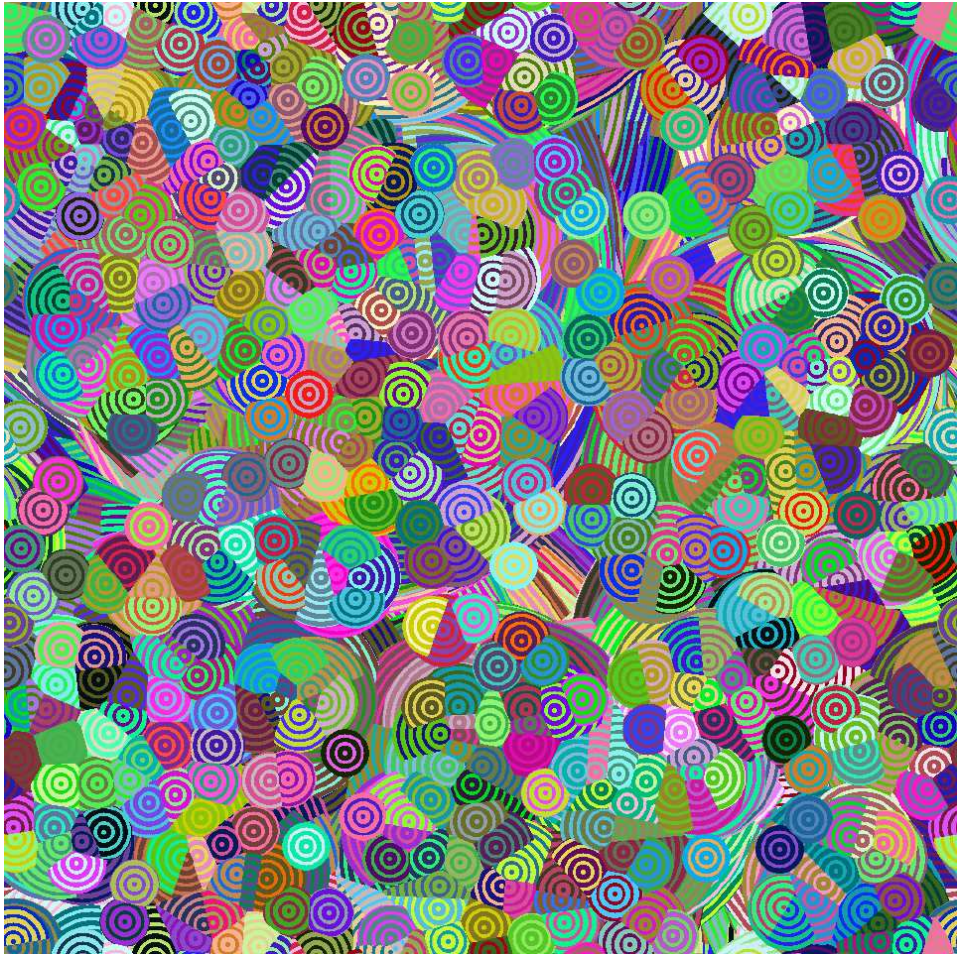

Fig. 1. *A stable allocation in the critical case $\alpha = 1$. (If you are looking at greyscale images, color versions are available at [www.math.ubc.ca/˜holroyd](www.math.ubc.ca/˜holroyd) or math.PR/0505668.)*

THEOREM 2 (Uniqueness for recurrent configurations). *Given discrete $\Xi \subseteq \mathbb{R}^d$, let $G_b$ be the graph with vertex set $\Xi$ and with an edge between every pair of centers at distance at most $b$ from each other. Suppose $\Xi$ is such that, for all $b > 0$, the graph $G_b$ is recurrent for simple symmetric random walk. Then for any $\alpha \in (0, \infty)$, there exists an $\mathcal{L}$-a.e. unique stable allocation.*

Now let $\Pi$ be a translation-invariant simple point process on $\mathbb{R}^d$, with intensity (i.e., mean points per unit volume) $\lambda \in (0, \infty)$ and law **P**. The *support* of $\Pi$ is the random set $[\Pi] = \{z \in \mathbb{R}^d : \Pi(\{z\}) = 1\}$. We consider stable allocations of the random set of centers $\Xi = [\Pi]$.

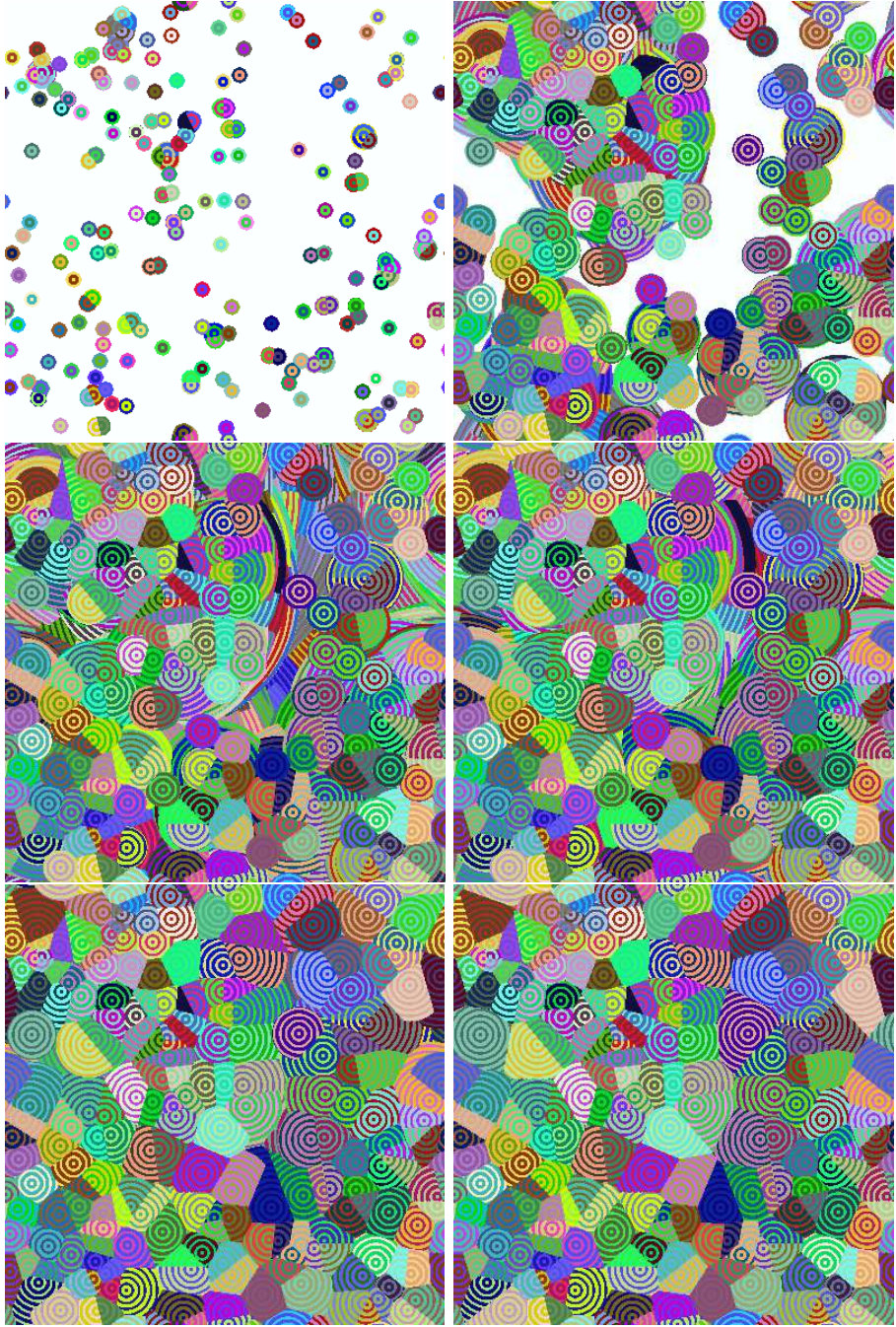

FIG. 2. *Stable allocations with appetites $\alpha = 0.2, 0.8, 1, 1.2, 2, \infty$. The centers are chosen uniformly at random in a 2-torus, with one center per unit area. Each territory is represented by concentric annuli in two colors.*

THEOREM 3 (Almost sure uniqueness). *For a translation-invariant point process $\Pi$ of finite intensity on $\mathbb{R}^d$ and for any $\alpha \in (0, \infty)$, there exists $\mathbf{P}$-a.s. an $\mathcal{L}$-a.e. unique stable allocation $\Psi_\Pi$ to $\Xi = [\Pi]$. Furthermore, $\Psi_\Pi$ can be chosen to be a measurable, isometry-equivariant factor of $\Pi$.*

Note that uniqueness can fail for a general set of centers, at least in the setting of sufficiently general spaces. For example, let $d = 1$, and $\Xi = \{\sum_{i=1}^n 1/i : n = 1, 3, 5, \ldots\}$, and consider instead of $\mathcal{L}$ the measure with a unit point mass at each of $\sum_{i=1}^n 1/i : n = 2, 4, 6, \ldots$. With $\alpha = 1$, there are multiple stable allocations: one assigns each mass of sites to the center immediately on its left; another assigns them to the right and leaves the center 1 unsated. We do not know whether uniqueness can fail for stable allocations of Lebesgue measure to a set of centers in $\mathbb{R}^d$.

Henceforth, we let $\Psi = \Psi_\Pi$ be as in Theorem 3.

THEOREM 4 (Phases). *Consider a stable allocation with appetite $\alpha$ to a point process of intensity $\lambda \in (0, \infty)$ which is ergodic under translations.*

(i) *If $\lambda\alpha < 1$ (subcritical), then a.s. all centers are sated but there is an infinite volume of unclaimed sites.*

(ii) *If $\lambda\alpha = 1$ (critical), then a.s. all centers are sated and $\mathcal{L}$-a.a. sites are claimed.*

(iii) *If $\lambda\alpha > 1$ (supercritical), then a.s. not all centers are sated but $\mathcal{L}$-a.a. sites are claimed.*

Note that by scaling $\mathbb{R}^d$ there is no loss of generality in assuming $\lambda = 1$ or $\alpha = 1$. [The assumption of ergodicity can of course be weakened to cover any mixture of ergodic point processes, all of whose intensities satisfy the appropriate condition in (i)–(iii).]

Note that the assertions of Theorem 4 may fail when generalized in the most obvious way to a general set of centers $\Xi$. Indeed, let $d = 1$ and $\alpha = 1$ and let $\Xi = \mathbb{Z} \setminus \{0\}$. The set $\Xi$ has "intensity" 1 in the sense that $\#\{\Xi \cap [-r, r]\} \sim \mathcal{L}[-r, r]$ as $r \to \infty$, so we might expect the configuration to be "critical." However, the stable allocation has a.e. site in $(-1/2, 1/2)$ unclaimed. (This may be deduced from a monotonicity result, Proposition 21, by comparison with the case $\Xi = \mathbb{Z}$.)

We now consider the geometry of territories.

THEOREM 5 (Geometry of territories). *Consider a stable allocation to a translation-invariant point process of finite intensity.*

(i) *All territories are bounded a.s.*

(ii) *Any bounded subset of $\mathbb{R}^d$ intersects only finitely many territories a.s.*

(iii) $\Psi_\Pi$ *may be chosen to be a stable isometry-equivariant factor of $\Pi$ in such a way that a.s. each territory is a union of finitely many open connected sets.*

We now restrict to the case when $\Pi$ is a homogeneous Poisson process with intensity 1. The following result contrasts with Theorem 5.

THEOREM 6 (Infinite demand). *Let $\Pi$ be a Poisson process and consider the critical case $\lambda = \alpha = 1$.*

  (i) *A.s. every center is desired by an infinite volume of sites.*
  (ii) *A.s. every site is coveted by infinitely many centers.*

[Note that (i) above holds trivially in the subcritical case, while (ii) holds trivially in the supercritical case, by Theorem 4.]

We now consider quantitative properties of stable allocations. The quantitative behavior turns out to be very different in the critical and noncritical phases. Broadly speaking, we find exponential tails in noncritical cases and power law tails in critical cases. A natural quantity to consider is

$$X = |\Psi_\Pi(0)|$$

(i.e., the distance from the origin to its center), where we take $|\infty| = \infty$. Let $\mathbf{E}$ be the expectation operator associated with $\mathbf{P}$.

THEOREM 7 (Critical lower bounds). *Let $\Pi$ be a Poisson process with intensity $\lambda = 1$.*

  (i) *For $d = 1, 2$ and $\alpha = 1$, we have $\mathbf{E}X^{d/2} = \infty$.*
  (ii) *For $d \geq 3$ and $\alpha = 1$, we have $\mathbf{E}X^d = \infty$.*

Proofs of the following two theorems appear in a second article by the same authors [8].

THEOREM II-1 (Critical upper bound). *Let $\Pi$ be a Poisson process with intensity $\lambda = 1$. For $d = 1$ and $\alpha = 1$, we have $\mathbf{E}X^{1/18} < \infty$.*

We have not been able to prove *any* quantitative upper bound for $X$ when $d \geq 2$.

THEOREM II-2 (Noncritical upper bounds). *Let $\Pi$ be a Poisson process with intensity $\lambda = 1$.*

  (i) *For all $d$ and $\alpha > 1$, we have $\mathbf{E}e^{cX^d} < \infty$;*
  (ii) *for all $d$ and $\alpha < 1$, we have $\mathbf{E}(e^{cX^d}; X < \infty) < \infty$,*

*where* $c = c(d, \alpha) > 0$.

It is straightforward to see that the bounds in Theorem II-2 are tight up to the value of $c$, since the distance from the origin to the *closest* center has tail decaying exponentially in volume. Some of the above bounds may be strengthened to bounds involving the radius of a typical territory; see [8] for more details.

REMARKS. (i) The concept of stability was introduced by Gale and Shapley [6] in the context of discrete "marriage problems," and has since been studied extensively; see [12] for more information. The model considered here appears to be the first extension to a continuum setting.

(ii) The cases $\alpha = 0, \infty$ are relatively simple. When $\alpha = 0$, every center's territory is either empty or consists only of the center itself. When $\alpha = \infty$, every center covets every site, so each $\psi(x)$ equals the closest center to $x$ and we obtain the Voronoi tessellation (see [14] for more information). Thus, we restrict our attention to the cases $\alpha \in (0, \infty)$.

(iii) A stable allocation has the following informal interpretation. Denote the ball $B(x, r) = \{y \in \mathbb{R}^d : |x - y| < r\}$. Fix $\Xi$ and $\alpha$, and imagine a sphere centered at each center, simultaneously growing linearly, so that, at time $t$, all spheres have radius $t$. Now let each center capture all sites encountered by its sphere, provided the center is not yet sated, and provided no other center captured the site first. It is easy to convince oneself that if this informal model exits, it should give rise to a stable allocation, and the territory captured by $\xi$ up to time $t$ should be $\psi^{-1}(\xi) \cap B(\xi, t)$. However, it is by no means clear a priori that the model does exist, since the proposed evolution involves long-range dependence. (See [1] for an example of a seemingly natural model which turns out not to exist in certain settings.) One consequence of proving existence and uniqueness of stable allocations is that the informal picture described does indeed make sense.

(iv) Our proof of existence (Theorem 1) will be a generalization of a proof in [10], which in turn is an adaptation of the algorithm introduced in [6]. In [10], stable allocations were used to construct certain shift-couplings of point processes.

(v) It should be noted that stable marriage problems do not in general have unique solutions; see [6]. In this case the key to uniqueness (Theorems 2 and 3) is that the preferences of sites and centers are "consistent," both being based on Euclidean distance.

(vi) The condition on $\Xi$ in Theorem 2 is natural when $d \leq 2$. In particular, by Corollary (3.52) in [16], it holds provided, for any $b < \infty$, we have

$$\#\{\{\xi, \xi'\} \subseteq \Xi \cap B(0, r) : |\xi - \xi'| < b\} = O(r^2 \log r) \qquad \text{as } r \to \infty.$$

(vii) In the course of proving our main theorems, we will obtain several other results which are interesting in their own right, and which may be summarized as follows. Stable allocations are monotonic in both $\Xi$ and $\alpha$ in the sense that adding extra centers or increasing the appetite makes sites happier but centers more unhappy (Propositions 21 and 22). Although stable allocations are unique only almost everywhere, there is a unique "canonical version" having the properties in Theorem 5 and other good properties (Theorem 24). In the second article [8] we prove that stable allocations are robust to far-away changes in $\Xi$, in the sense that if a sequence of sets of centers converges weakly in the discrete topology, then the sequence of stable allocations converges almost everywhere. The supercritical phase is rigid in the sense that if $\alpha$ is only just greater than one, a large proportion of centers are sated.

(viii) Note that there is a partial symmetry between the roles of sites and centers in the definition of stability. In light of this, it is not surprising that many of our results (and the arguments used to prove them) come in "dual" pairs obtained by exchanging the roles played by sites and centers.

**2. Existence.** In order to prove Theorem 1, we shall construct an allocation $\psi$ using a variant of the Gale–Shapley stable marriage algorithm [6]. The algorithm has two different versions, called "site-optimal" and "center-optimal." (The reason for the names will become apparent in Lemma 15.) For our present purpose, one version will be sufficient, but later in the proof of Theorem 3 we shall have occasion to consider both.

*Site-optimal Gale–Shapley algorithm.* Let $W$ be the set of all sites which are equidistant from two or more centers. Since $\Xi$ is countable, $W$ is $\mathcal{L}$-null. For convenience, we take $\psi(x) = \Delta$ for all $x \in W$.

We now construct $\psi$ by means of a sequence of stages. For each positive integer $n$, *stage* $n$ consists of two parts as follows:

(a) Each site $x \notin W$ *applies* to the closest center to $x$ which has not rejected $x$ at any earlier stage.

(b) For each center $\xi$, let $A_n(\xi)$ be the set of sites which applied to $\xi$ in stage $n$ (a), and define the *rejection radius*

$$r_n(\xi) = \inf\{r : \mathcal{L}(A_n(\xi) \cap B(\xi, r)) \geq \alpha\},$$

where the infimum of the empty set is taken to be $\infty$. Then $\xi$ *shortlists* all sites in $A_n(\xi) \cap B(\xi, r_n(\xi))$, and *rejects* all sites in $A_n(\xi) \setminus B(\xi, r_n(\xi))$.

We now describe $\psi$. Consider a site $x \notin W$. Since $\Xi$ is discrete, the following is clear. Either $x$ is rejected by every center (in order of increasing distance from $x$) or for some center $\xi$ and some stage $n$, $x$ is shortlisted by $\xi$ at all stages after $n$. In the former case we put $\psi(x) = \infty$ (so $x$ is unclaimed); in the latter case we put $\psi(x) = \xi$.

PROOF OF THEOREM 1. Let $\psi$ be constructed according to the site-optimal Gale–Shapley algorithm above. We first check that $\psi$ is an allocation. Let $S_n(\xi)$ be the set of sites shortlisted by a center $\xi$ at stage $n$. By the construction in (b) and the intermediate value theorem, we have $\mathcal{L}(S_n(\xi)) \leq \alpha$. But by the definition of $\psi$ above, we have $\psi^{-1}(\xi) = \limsup_{n\to\infty} S_n(\xi) = \liminf_{n\to\infty} S_n(\xi)$, so, by Fatou's lemma, we have $\mathcal{L}(\psi^{-1}(\xi)) \leq \alpha$, as required. Note also that if a center $\xi$ ever rejects sites (at stage $n$ say), then we must have $\mathcal{L}(S_m(\xi)) = \alpha$ for all stages $m \geq n$. Hence, an unsated center never rejected any sites.

We must check that $\psi$ is stable. Consider a site $x$ and a center $\xi$.

Suppose $x$ desires $\xi$. If $x$ is unclaimed, then it was rejected by all centers. On the other hand, if $x$ is claimed and $|x - \xi| < |x - \psi(x)|$, then the site $x$ applied to $\psi(x)$ at some stage, therefore, it was rejected by the closer center $\xi$ at some earlier stage. In either case $\xi$ rejected $x$.

Now suppose $\xi$ covets $x$. If $|x - \xi| < |x' - \xi|$ for some $x' \in \psi^{-1}(\xi)$, then $\xi$ never rejected $x'$, so $x'$ was never outside $\xi$'s rejection radius, and the same holds for the closer site $x$. On the other hand, if $\xi$ is unsated, then it never rejected any sites. In either case $\xi$ did not reject $x$.

We deduce as required that $(x, \xi)$ cannot be an unstable pair. $\square$

**3. Uniqueness for recurrent configurations.** In this section we prove Theorem 2. We therefore consider a set $\Xi$ satisfying the assumption of the theorem. Without loss of generality (by scaling $\mathbb{R}^d$), we may also assume that $\alpha = 1$. We work by contradiction. Suppose that $\psi_1$ and $\psi_2$ are two stable allocations to $\Xi$ which differ on an $\mathcal{L}$-nonnull set. We will use them to define a certain flow on a directed graph with vertex set $\Xi$, and then show that this flow has properties which are impossible when $\Xi$ satisfies the given condition.

For each center $\xi$, we start by defining several sets. Let the set of *disputed* sites for $\xi$ be the symmetric set difference

$$D(\xi) = [\psi_1^{-1}(\xi) \triangle \psi_2^{-1}(\xi)] \setminus [\psi_1^{-1}(\Delta) \cup \psi_2^{-1}(\Delta)].$$

Let

$$N(\xi) = D(\xi) \cap B(\xi, r(\xi)),$$

where

$$r(\xi) = \sup\{r : \mathcal{L}[(\psi_1^{-1}(\xi) \cup \psi_2^{-1}(\xi)) \cap B(\xi, r)] < 1\}.$$

This is the set of disputed sites that $\xi$ prefers. Also define

$$M(\xi) = D(\xi) \cap \{x : |x - \xi| < \max[|x - \psi_1(x)|, |x - \psi_2(x)|]\}.$$

This is the set of disputed sites that prefer $\xi$.

Define for every pair of centers $\xi, \xi'$,

$$F(\xi, \xi') = N(\xi) \cap M(\xi').$$

Let $G$ be the directed graph with vertex set $\Xi$ and a directed edge from $\xi$ to $\xi'$ whenever $\mathcal{L}F(\xi, \xi') > 0$. We label the directed edge $(\xi, \xi')$ with the quantity $\mathcal{L}F(\xi, \xi')$. We will regard these labels as capacities. We first show that $G$ has no directed cycles.

LEMMA 8. *For each center $\xi$,*

$$\mathcal{L}F(\xi, \xi) = 0.$$

*Thus, $G$ has no cycles of length one (self-loops).*

PROOF. Suppose that $\mathcal{L}F(\xi, \xi) > 0$. Then, without loss of generality, there exist a region $A$ of positive volume and a center $\xi' \neq \xi$ such that:

(i) $A \subseteq N(\xi) \cap M(\xi)$,
(ii) $A \subseteq \psi_1^{-1}(\xi)$ and
(iii) $A \subseteq \psi_2^{-1}(\xi')$.

We claim that either

(a) $\mathcal{L}(\psi_2^{-1}(\xi) \setminus B(\xi, r(\xi))) \geq \mathcal{L}A \ (> 0)$, or
(b) $\mathcal{L}(\psi_2^{-1}(\xi)) < 1$.

To see this, suppose, contrary to (b), that $\mathcal{L}\psi_2^{-1}(\xi) = 1$. The definition of $r(\xi)$ then implies

$$1 = \mathcal{L}([\psi_1^{-1}(\xi) \cup \psi_2^{-1}(\xi)] \cap B(\xi, r(\xi)))$$
$$= \mathcal{L}(\psi_2^{-1}(\xi) \cap B(\xi, r(\xi))) + \mathcal{L}([\psi_1^{-1}(\xi) \setminus \psi_2^{-1}(\xi)] \cap B(\xi, r(\xi))).$$

Conditions (i), (ii) and (iii) above imply $A \subseteq [\psi_1^{-1}(\xi) \setminus \psi_2^{-1}(\xi)] \cap B(\xi, r(\xi))$, therefore, $\mathcal{L}[\psi_2^{-1}(\xi) \cap B(\xi, r(\xi))] \leq 1 - \mathcal{L}A$. Since $\mathcal{L}\psi_2^{-1}(\xi) = 1$, this implies (a) above.

Now consider a site $x \in A$. We have that $x \in N(\xi) \subseteq B(\xi, r(\xi))$. In case (a), for any $y \in \psi_2^{-1}(\xi) \setminus B(\xi, r(\xi))$, we have

$$|x - \xi| < |y - \xi|.$$

Hence, in either case (a) or (b), $\xi$ covets $x$ under $\psi_2$. Also, since $x \in M(\xi)$, we have that

$$|x - \xi| < |x - \xi'|,$$

therefore, $x$ desires $\xi$ under $\psi_2$. Thus, $(x, \xi)$ is an unstable pair for $\psi_2$, a contradiction. □

LEMMA 9 (Monotonicity of distances). *For centers $\xi_1 \neq \xi_2$ and $x_1 \in F(\xi_1, \xi_2)$, we have the following:*

(i) $|\xi_1 - x_1| > |x_1 - \xi_2|$, *and*
(ii) *for any center $\xi_0$ and $x_0 \in F(\xi_0, \xi_1) \setminus N(\xi_1)$,*

$$|x_0 - \xi_1| > |\xi_1 - x_1|.$$

[Note that, by Lemma 8, the condition in (ii) above applies to a.e. site $x_0 \in F(\xi_0, \xi_1)$.]

PROOF. Since $x_1 \in M(\xi_2)$, statement (i) follows from the definition of $M(\xi_2)$. For (ii), since $x_1 \in N(\xi_1)$ and $x_0 \in D(\xi_1) \setminus N(\xi_1)$, the definition of $N(\xi_1)$ implies that $|x_1 - \xi_1| < r(\xi_1) \leq |x_0 - \xi_1|$. □

LEMMA 10. *G has no directed cycles.*

PROOF. By Lemma 8, there are no cycles of length 1. Suppose there is a cycle with vertices $\xi_0, \xi_1, \xi_2, \ldots, \xi_n = \xi_0$. Then for all $i \in \{0, \ldots, n-1\}$, there exists $x_i \in F(\xi_i, \xi_{i+1})$. By Lemma 8, we can further require that $x_i \notin N(\xi_{i+1})$. Then Lemma 9 yields

$$|\xi_0 - x_0| > |x_0 - \xi_1| > |\xi_1 - x_1| > \cdots > |x_{n-1} - \xi_n| > |\xi_n - x_0|$$
$$= |\xi_0 - x_0|. \qquad \square$$

LEMMA 11 (Increasing capacity). *For any center $\xi$, we have*

$$\sum_{\xi' \in \Xi} \mathcal{L}F(\xi', \xi) \leq \sum_{\xi'' \in \Xi} \mathcal{L}F(\xi, \xi'').$$

PROOF. By the definitions of $F(\xi, \xi'), N(\xi)$ and $M(\xi)$, we have that

(1) $$\bigsqcup_{\xi' \in \Xi} F(\xi', \xi) \subseteq M(\xi)$$

and

(2) $$\bigsqcup_{\xi'' \in \Xi} F(\xi, \xi'') = N(\xi).$$

(Here ⊔ denotes disjoint union.) We claim now that

$$\mathcal{L}M(\xi) \leq \mathcal{L}N(\xi),$$

which combined with (1) and (2) proves the lemma. To prove the claim, we consider two cases. First, if

$$\mathcal{L}N(\xi) < 1 - \mathcal{L}(\psi_1^{-1}(\xi) \cap \psi_2^{-1}(\xi)),$$

then $r(\xi) = \infty$, so $N(\xi) = D(\xi)$, and by Lemma 8, it follows that $\mathcal{L}M(\xi) = 0$, establishing the claim. Second, suppose

$$(3) \qquad \mathcal{L}N(\xi) = 1 - \mathcal{L}(\psi_1^{-1}(\xi) \cap \psi_2^{-1}(\xi)).$$

By Lemma 8, we have

$$2\mathcal{L}(\psi_1^{-1}(\xi) \cap \psi_2^{-1}(\xi)) + \mathcal{L}N(\xi) + \mathcal{L}M(\xi) \leq \mathcal{L}\psi_1^{-1}(\xi) + \mathcal{L}\psi_2^{-1}(\xi) \leq 2,$$

and combining this with (3) yields

$$\mathcal{L}M(\xi) \leq 1 - \mathcal{L}(\psi_1^{-1}(\xi) \cap \psi_2^{-1}(\xi)),$$

again establishing the claim. □

Supposing that $\psi_1$ and $\psi_2$ do not agree $\mathcal{L}$-a.e., there exists some center $\zeta$ with $\mathcal{L}D(\zeta) > 0$, and therefore, $\mathcal{L}N(\zeta) > 0$. Any site $x \in D(\zeta)$ lies in $M(\zeta')$ for some center $\zeta'$. Therefore, $G$ has some edge $(\zeta, \zeta')$. Now fix a directed edge $(\zeta, \zeta')$ of $G$ and let $H$ be the directed subgraph consisting of this edge together with all vertices and edges of $G$ that can be reached by directed paths from $\zeta'$.

LEMMA 12 (Bounded edges and summable paths). *Let $H$ be a graph as defined above.*

  (i) *The edges of $H$ have uniformly bounded length.*
  (ii) *There exists $K < \infty$ such that, for any directed path $(\xi_0, \xi_1, \xi_2, \ldots)$ in $H$, we have*

$$\sum_{i=0}^{\infty} \mathcal{L}F(\xi_i, \xi_{i+1}) \leq K.$$

PROOF. We start by proving (ii). Let $(\xi_0, \xi_1, \xi_2, \ldots)$ be a directed path. Without loss of generality, we may assume that $\xi_0 = \zeta$ and $\xi_1 = \zeta'$. Define

$$a_i = \operatorname*{ess\,inf}_{z \in F(\xi_i, \xi_{i+1})} |\xi_i - z|$$

and

$$b_i = \operatorname*{ess\,sup}_{z \in F(\xi_i, \xi_{i+1})} |\xi_i - z|.$$

Clearly, $b_i \geq a_i$ for each $i$, and we claim that also

$$(4) \qquad a_i \geq b_{i+1}, \qquad i = 0, 1, \ldots.$$

To prove this, note that, by Lemma 9 (and Lemma 8), for $\mathcal{L}$-a.e. $x_i \in F(\xi_i, \xi_{i+1})$ and every $x_{i+1} \in F(\xi_{i+1}, \xi_{i+2})$, we have

$$(5) \qquad |\xi_i - x_i| \geq |x_i - \xi_{i+1}| \geq |\xi_{i+1} - x_{i+1}|;$$

then use the definitions of $a_i, b_{i+1}$.

Let $c = c(d) < \infty$ be the volume of the unit ball. Using (4), we have

$$\sum_{i=0}^{\infty} \mathcal{L}F(\xi_i, \xi_{i+1}) \leq \mathcal{L}D(\xi_0) + \sum_{i=1}^{\infty}(cb_i^d - ca_i^d)$$

$$\leq 2 + c\sum_{i=1}^{\infty}(a_{i-1}^d - a_i^d)$$

$$= 2 + ca_0^d.$$

But $a_0 = \mathrm{ess\,inf}_{z \in F(\zeta,\zeta')}|\zeta - z|$, and this quantity is finite and independent of the choice of path, so we obtain (ii).

Turning to (i), note that for a directed path as above, for $i \geq 1$ and $\mathcal{L}$-a.e. $x_i \in F(\xi_i, \xi_{i+1})$, we have from (5) and (4) that

$$|\xi_i - \xi_{i+1}| \leq |\xi_i - x_i| + |x_i - \xi_{i+1}|$$

$$\leq 2|\xi_i - x_i|$$

$$\leq 2b_i \leq 2a_0.$$

Since every edge of $H$ lies in some directed path starting with the edge $(\zeta, \zeta')$, (i) follows. □

We will now deduce Theorem 2 from the following result on random walks. The proof of Proposition 13 uses standard technology, and is deferred until later in this section.

PROPOSITION 13 (Transience from capacities).  *Let $H = (V, E)$ be a locally finite, directed graph with no directed cycles, and with a distinguished vertex $\rho \in V$. Suppose there exists a "capacity" function $q : E \to [0, \infty)$ which satisfies the following three conditions:*

(i) *For every vertex $v \neq \rho$, we have*

$$\sum_{u:(u,v) \in E} q(u,v) \leq \sum_{w:(v,w) \in E} q(v,w).$$

(ii) $\sum_{u:(u,\rho) \in E} q(u,\rho) = 0 \quad \text{and} \quad \sum_{w:(\rho,w) \in E} q(\rho, w) > 0.$

(iii) *There exists $K < \infty$ such that, for every directed path $(u_0, u_1, \ldots)$ in $H$, starting at $u_0 = \rho$, we have*

$$\sum_{i=0}^{\infty} q(u_i, u_{i+1}) \leq K.$$

*Then simple symmetric random walk on the undirected version of $H$ starting at $\rho$ is transient.*

PROOF OF THEOREM 2.  Suppose $\Xi$ has two stable allocations which differ on an $\mathcal{L}$-nonnull set, and construct the graph $H$ as described before Lemma 12 above. By the construction of $H$ and Lemmas 10, 11 and 12(ii), $H$ satisfies the assumptions of Proposition 13 with $q(\xi, \xi') = \mathcal{L}F(\xi, \xi')$, hence, it is transient (for simple symmetric random walk on the undirected version of the graph). However, by Lemma 12(i), the assumption on $\Xi$ and the fact that a subgraph of a recurrent graph is recurrent (Corollary (3.48) of [16]), $H$ is recurrent. This is a contradiction. □

PROOF OF PROPOSITION 13.  Without loss of generality, we may assume $q$ is normalized so that $\sum_{w:(\rho,w)\in E} q(\rho, w) = 1$. By a *unit flow* on $H$, we mean a function $f : E \to [0, \infty)$ satisfying the following:

(i) for every vertex $v \neq \rho$, we have
$$\sum_{u:(u,v)\in E} f(u,v) = \sum_{w:(v,w)\in E} f(v,w);$$

(ii) $\sum_{u:(u,\rho)\in E} f(u,\rho) = 0$ and $\sum_{w:(\rho,w)\in E} f(\rho,w) = 1.$

First we claim that there exists a unit flow $f$ with $f \leq q$ on $E$. This is proved using the max-flow/min-cut theorem (see, e.g., [15], Theorem 2.19). As a notational convenience we extend $q$ to be 0 on $(V \times V) \setminus E$. For any finite set of vertices $S \ni \rho$, conditions (ii) and (iii) of the proposition together with the normalization assumption yield
$$\sum_{v\in S}\sum_{w\notin S} q(v,w) \geq \sum_{v\in S}\sum_{w\in V} [q(v,w) - q(w,v)] \geq 1,$$

so the claim follows. Now from condition (iv) of the proposition, we have also

(6) $$\sum_{i=0}^{\infty} f(u_i, u_{i+1}) \leq K$$

for every directed path $(u_0, u_1, \ldots)$ starting at $\rho$.

If $(u_0, u_1, \ldots)$ is a directed path starting at $\rho$, then there is a unit flow which takes the value 1 on each edge of the path and 0 elsewhere. Let $\Gamma$ be the set of all unit flows of this form. By Proposition 2.20 of [15], for any acyclic unit flow $f$, there exists a unit flow $f_0$ satisfying $f_0 \leq f$ on $E$, and a probability measure $\mu$ on $\Gamma$, such that, for all $e \in E$,

(7) $$f_0(e) = \int_\Gamma \gamma(e)\, d\mu(\gamma).$$

Now, apply this to the unit flow $f$ constructed earlier, which is acyclic since $f \leq q$. The summability condition on paths (6) implies that, for any $\gamma \in \Gamma$, we have
$$\sum_{e \in E} f_0(e)\gamma(e) \leq K.$$
Therefore, by (7) and Fubini's theorem,
$$\sum_{e \in E} f_0(e)^2 = \sum_{e \in E} f_0(e) \int_\Gamma \gamma(e)\, d\mu(\gamma) \leq \int K\, d\mu = K.$$
So $f_0$ is a unit flow of finite energy, and it follows (see, e.g., [4]) that $H$ is transient as required. □

**4. Almost sure uniqueness.** In order to prove Theorem 3, we introduce a second version of the Gale–Shapley algorithm which constructs a potentially different allocation. Our approach to proving uniqueness will be to show that the two Gale–Shapley allocations are extremal (in an appropriate sense), and then to show that they in fact coincide.

*Center-optimal Gale–Shapley algorithm.* Let $W$ be the set of all sites which are equidistant from two or more centers, and take $\psi(x) = \Delta$ for all $x \in W$. For each positive integer $n$, *stage $n$* consists of two parts as follows:

(a) For each center $\xi$, let $R_n(\xi)$ be the set of sites which have rejected $\xi$ at some earlier stage, and define the *application radius*
$$a_n(\xi) = \inf\{a : \mathcal{L}(B(\xi, a) \setminus R_n(\xi)) \geq \alpha\}.$$
The center $\xi$ *applies* to every site $x \in B(\xi, a_n(\xi)) \setminus R_n(\xi) \setminus W$.

(b) Each site $x \notin W$ *shortlists* the closest center of those which applied to $x$ in stage $n$ (a) (if any), and *rejects* any others which applied.

We now describe $\psi$. Consider a site $x \notin W$. If no center ever applies to $x$, we put $\psi(x) = \infty$. Otherwise, if $x$ shortlists some center at some stage, then $x$ will shortlist this or a closer center at all later stages. Since $\Xi$ is discrete, it follows that, for some center $\xi$ and some stage $n$, $x$ shortlists $\xi$ at all stages after $n$. In this case we put $\psi(x) = \xi$.

LEMMA 14. *For any $\Xi$, the center-optimal Gale–Shapley algorithm yields a stable allocation.*

PROOF. We mimic the proof of Theorem 1 in the site-optimal case. To see that $\psi$ is an allocation, note that the set $A_n(\xi)$ of sites which $\xi$ applies to at stage $n$ satisfies $\mathcal{L}(A_n) \leq \alpha$, and $\psi^{-1}(\xi) = \liminf_{n \to \infty} A_n(\xi) = \limsup_{n \to \infty} A_n(\xi)$, and apply Fatou's lemma.

To check stability, consider a pair $(x, \xi)$. If $x$ desires $\xi$, then $\xi$ did not apply to $x$, while if $\xi$ covets $x$, then $\xi$ did apply to $x$. Hence, $(x, \xi)$ is not unstable. $\square$

For any allocation $\psi$ and any $r \in [0, \infty]$, define the following quantities. For any site $x \notin \psi^{-1}(\Delta)$,

$$g_x(\psi, r) = \mathbf{1}[|x - \psi(x)| < r];$$

and for any center $\xi$,

$$\gamma_\xi(\psi, r) = \mathcal{L}[\psi^{-1}(\xi) \cap B(\xi, r)].$$

As $r$ varies, the above quantities measure how good the allocation is from the point of view of a particular site or center respectively.

The following result states that the site-optimal Gale–Shapley allocation is the best for sites and the worst for centers, and vice versa for the center-optimal allocation. (The analogous fact for discrete stable marriage problems is well known: the algorithm in which boys propose gives the best possible stable system of marriages for boys and the worst possible for girls—see [6, 12].)

LEMMA 15 (Optimality). *Fix a set of centers $\Xi$, let $\psi$ be any stable allocation, and let $\psi^{\mathrm{GS}}, \psi^{\mathrm{GSc}}$ be respectively the site-optimal and center-optimal Gale–Shapley allocations to the same $\Xi$. For all $r \in [0, \infty]$, we have the following:*

(i) *for $\mathcal{L}$-a.e. site $x$,*

$$g_x(\psi^{\mathrm{GS}}, r) \geq g_x(\psi, r) \geq g_x(\psi^{\mathrm{GSc}}, r);$$

(ii) *for every center $\xi$,*

$$\gamma_\xi(\psi^{\mathrm{GSc}}, r) \geq \gamma_\xi(\psi, r) \geq \gamma_\xi(\psi^{\mathrm{GS}}, r).$$

The proof of Lemma 15 is deferred until the end of this section.

*Stable allocations to a point process.* Now let $\Pi$ be a translation-invariant simple point process on $\mathbb{R}^d$ with intensity $\lambda \in (0, \infty)$, law $\mathbf{P}$ and expectation operator $\mathbf{E}$. The assumption of finite intensity ensures that $[\Pi]$ is discrete a.s. A *factor allocation* is a measurable map $\Psi_\bullet$ which assigns to $\mathbf{P}$-a.e. point measure $\pi$ an allocation $\Psi_\pi$ of $[\pi]$, and which is translation-equivariant in the sense that, for any translation $T$ of $\mathbb{R}^d$, if $\Psi_\pi(x) = \xi$, then $\Psi_{T\pi}(Tx) = T\xi$ (here it is understood that $\Psi_{T\pi}$ is defined if and only if $\Psi_\pi$ is). A factor allocation is *isometry-equivariant* if the above holds for all isometries $T$. Note that the two Gale–Shapley algorithms clearly yield isometry-equivariant factor allocations; we denote them $\Psi_\bullet^{\mathrm{GS}}, \Psi_\bullet^{\mathrm{GSc}}$, respectively. For a factor allocation $\Psi_\bullet$, we shall consider the random allocation $\Psi_\Pi$, and when there is no risk of confusion, we denote it simply $\Psi$. Note that the joint law of the pair $(\Pi, \Psi_\Pi)$ is then translation-invariant.

Let $\Pi^*$ be the Palm version of the point process $\Pi$, with law $\mathbf{P}^*$ and expectation operator $\mathbf{E}^*$. (Recall that $\Pi^*$ may be thought of as $\Pi$ conditioned to have a center at the origin, and that when $\Pi$ is a Poisson process, $\Pi^*$ is a Poisson process with an added center at the origin; see [11], Chapter 11 for details.) By results of Thorisson [17] (see also [10] and [11]), $\Pi^*$ and $\Pi$ may be coupled so that one is almost surely a translation of the other. Hence, if $\Psi_\bullet$ is a factor allocation, then $\Psi_{\Pi^*}$ is defined $\mathbf{P}^*$-a.s.

LEMMA 16 (Site-center consistency). *If $\Pi$ is a translation-invariant point process of intensity $\lambda \in (0, \infty)$ and $\Psi_\bullet$ is a factor allocation of $\Pi$, then*

$$\mathbf{E} g_0(\Psi_\Pi, r) = \lambda \mathbf{E}^* \gamma_0(\Psi_{\Pi^*}, r)$$

*for all $r$.*

The key to Lemma 16 is that the preferences of sites and centers expressed in the functions $g, \gamma$ are "consistent"—both being defined in terms of Euclidean distance, so that on average they coincide. The proof is deferred to the end of this section.

PROOF OF THEOREM 3. Note that $\Psi_\bullet^{\mathrm{GS}}, \Psi_\bullet^{\mathrm{GSc}}$ are clearly isometry-equivariant factor allocations. Taking expectations in Lemma 15 and applying Lemma 16, we obtain

$$\mathbf{E} g_0(\Psi_\Pi^{\mathrm{GS}}, r) \geq \mathbf{E} g_0(\Psi_\Pi^{\mathrm{GSc}}, r) = \lambda \mathbf{E}^* \gamma_0(\Psi_{\Pi^*}^{\mathrm{GSc}}, r)$$
$$\geq \lambda \mathbf{E}^* \gamma_0(\Psi_{\Pi^*}^{\mathrm{GS}}, r) = \mathbf{E} g_0(\Psi_\Pi^{\mathrm{GS}}, r).$$

Hence, we have equality throughout. Hence, by translation-invariance, we have

$$\mathbf{E} g_x(\Psi_\Pi^{\mathrm{GS}}, r) = \mathbf{E} g_x(\Psi_\Pi^{\mathrm{GSc}}, r)$$

for all $r$ and all sites $x$. Also, by Lemma 15, we have

$$g_x(\Psi_\Pi^{\mathrm{GS}}, r) \geq g_x(\Psi_\Pi^{\mathrm{GSc}}, r) \qquad \text{for all } r,$$

wherever both quantities are defined. Hence, for all sites $x$ we have $\mathbf{P}$-a.s.

(8) $\qquad g_x(\Psi_\Pi^{\mathrm{GS}}, r) = g_x(\Psi_\Pi^{\mathrm{GSc}}, r) \qquad \text{for all } r,$

wherever both quantities are defined. Therefore, by Fubini's theorem, $\mathbf{P}$-a.s., we have that (8) holds for $\mathcal{L}$-a.e. $x$.

Provided $x$ is not equidistant from two centers, (8) implies that $\Psi^{\mathrm{GS}}(x) = \Psi^{\mathrm{GSc}}(x)$. Hence, **P**-a.s., for $\mathcal{L}$-a.e. $x$, we have $\Psi^{\mathrm{GS}}(x) = \Psi^{\mathrm{GSc}}(x)$. Finally, using Lemma 15 again, **P**-a.s., we have, for any stable allocation $\psi$ of $\Xi = [\Pi]$ and $\mathcal{L}$-a.e. $x$,

$$g_x(\Psi_\Pi^{\mathrm{GS}}, r) \geq g_x(\psi, r) \geq g_x(\Psi_\Pi^{\mathrm{GSc}}, r) = g_x(\Psi_\Pi^{\mathrm{GS}}, r) \qquad \text{for all } r,$$

and so we have equality throughout. Thus, any two stable allocations to $[\Pi]$ agree for $\mathcal{L}$-a.e. $\square$

Finally, in this section we prove Lemmas 15 and 16. The proof of Lemma 16 proceeds via the following lemma (see [2] and [7] for generalizations).

LEMMA 17 (Mass-transport principle). *Let $m : \mathbb{Z}^d \times \mathbb{Z}^d \to [0, \infty]$ be a function satisfying $m(u+z, v+z) = m(u,v)$ for all $z \in \mathbb{Z}^d$. Then*

$$\sum_{v \in \mathbb{Z}^d} m(0, v) = \sum_{u \in \mathbb{Z}^d} m(u, 0).$$

PROOF.

$$\sum_{v \in \mathbb{Z}^d} m(0, v) = \sum_{v \in \mathbb{Z}^d} m(-v, 0) = \sum_{u \in \mathbb{Z}^d} m(u, 0). \qquad \square$$

For $z \in \mathbb{Z}^d$, we define the unit cube $Q_z = [0, 1)^d + z \subset \mathbb{R}^d$.

PROOF OF LEMMA 16. We shall apply the mass-transport principle (Lemma 17) to

$$m(u, v) = \mathbf{E}\mathcal{L}\{x \in Q_u : \Psi_\Pi(x) \in Q_v, |x - \Psi_\Pi(x)| < r\}.$$

Using Fubini's theorem and translation-invariance, we have

$$\sum_{v \in \mathbb{Z}^d} m(0, v) = \mathbf{E}\mathcal{L}\{x \in Q_0 : |x - \Psi_\Pi(x)| < r\} = \mathbf{P}(|0 - \Psi_\Pi(0)| < r),$$

while by a standard property of the Palm process,

$$\sum_{u \in \mathbb{Z}^d} m(u, 0) = \mathbf{E} \sum_{\xi \in [\Pi] \cap Q_0} \mathcal{L}[\Psi_\Pi^{-1}(\xi) \cap B(\xi, r)] = \lambda \mathbf{E}^* \mathcal{L}[\Psi_{\Pi^*}^{-1}(0) \cap B(0, r)].$$

So the result follows. $\square$

The proof of Lemma 15 depends on the two lemmas below. It is convenient to prove the first in a general form which will also be useful in proving monotonicity results in Section 6. For the present purpose, we may assume $\Xi_0 = \Xi$ and $\alpha_0 = \alpha$.

LEMMA 18. *Let $\Xi_0 \subseteq \Xi$ be two sets of centers, and let $\alpha_0 \geq \alpha$. Let $\psi$ be any stable allocation to $\Xi_0$ with appetite $\alpha_0$, and consider the site-optimal Gale–Shapley algorithm for $\Xi$ with appetite $\alpha$. Then for $\mathcal{L}$-a.e. site $x$, we have that $\psi(x)$ never rejects $x$.*

LEMMA 19. *Let $\Xi$ be a set of centers. Let $\psi$ be any stable allocation to $\Xi$ and consider the center-optimal Gale–Shapley algorithm for $\Xi$. Then for $\mathcal{L}$-a.e. site $x$, we have that $x$ never rejects $\psi(x)$.*

PROOF OF LEMMA 18. The proof is by induction on the stage. Suppose that the statement is false, and let $n$ be the first stage at which it is violated. So at stage $n$, some center $\xi \in \Xi_0$ rejects a set of sites $R \subseteq \psi^{-1}(\xi)$ with $\mathcal{L}R = \varepsilon > 0$. Then at stage $n$, the center must have shortlisted a set $S$ with $\mathcal{L}S = \alpha_0$, and with every site in $S$ closer to $\xi$ than some site in $R$ is to $\xi$. Since $\psi$ is an allocation, we have $\mathcal{L}\psi^{-1}(\xi) \leq \alpha \leq \alpha_0$, so there must exist a set $T \subseteq S$, disjoint from $\psi^{-1}(\xi)$, with $\mathcal{L}T \geq \varepsilon$. Since $\xi$ shortlisted $S$ at stage $n$, every site $x \in T$ must earlier have been rejected by all closer centers in $\Xi$. Hence, by the assumption, for $\mathcal{L}$-a.e. $x \in T$, we have $|\psi(x) - x| > |\xi - x|$. But this implies that $(x, \xi)$ is an unstable pair for $\psi$. □

PROOF OF LEMMA 19. The proof is by induction on the stage. Suppose that the statement is false, and let $n$ be the first stage at which it is violated. So at stage $n$, some nonnull set of sites $R \subseteq \psi^{-1}(\xi)$ all reject some center $\xi$. Therefore, each site in $R$ shortlisted some closer center at stage $n$, and we may therefore find some center $\xi'$ and some nonnull set $T \subseteq R$ such that every site in $T$ shortlisted $\xi'$. In particular, $\xi'$ applied to $T$ at stage $n$. Now consider $\psi^{-1}(\xi')$. If every site in $\psi^{-1}(\xi')$ is closer to $\xi'$ than $T$ is, then $\xi'$ must earlier have rejected some nonnull subset of $\psi^{-1}(\xi')$, which contradicts our assumption. Therefore, there exist sites $x \in R'$ and $y \in \psi^{-1}(\xi')$ with $|y - \xi'| > |x - \xi'|$. But then $(x, \xi')$ is an unstable pair for $\psi$. □

PROOF OF LEMMA 15. The first inequality in (i) is equivalent to the assertion that, for $\mathcal{L}$-a.e. $x$, we have

$$(9) \qquad |x - \psi^{\mathrm{GS}}(x)| \leq |x - \psi(x)|,$$

where we write $|x - \infty| = \infty$. Note that in the site-optimal Gale–Shapley algorithm, for any $x$ not equidistant from two centers, we have that $\psi^{\mathrm{GS}}(x)$ equals the closest center to $x$ which never rejects $x$, or $\infty$ if all centers reject $x$. Therefore, the first inequality in (i) follows from Lemma 18 with $\Xi_0 = \Xi$ and $\alpha_0 = \alpha$.

Next we consider the second inequality in (ii). Suppose for a contradiction that $\gamma_\xi(\psi, r) < \gamma_\xi(\psi^{\mathrm{GS}}, r)$. Since $\psi^{\mathrm{GS}}$ is an allocation, this implies that, in $\psi$, the center $\xi$ is unsated or has some site outside $B(\xi, r)$ in its territory. Also, there exists an $\mathcal{L}$-nonnull set $T \subseteq B(\xi, r) \cap [(\psi^{\mathrm{GS}})^{-1}(\xi) \setminus \psi^{-1}(\xi)]$. From (9), for $\mathcal{L}$-a.e. site $x \in T$, we have $|x - \psi^{\mathrm{GS}}(x)| < |x - \xi|$. But now $(x, \xi)$ is unstable for $\psi$.

Next we turn to the first inequality in (ii). Suppose, on the contrary, that $\gamma_\xi(\psi^{\mathrm{GSc}}, r) < \gamma_\xi(\psi, r)$. Since $\psi^{\mathrm{GSc}}$ is an allocation, this implies that, in $\psi^{\mathrm{GSc}}$, the center $\xi$ is unsated or has some site outside $B(\xi, r)$ in its territory. Also, there exists an $\mathcal{L}$-nonnull set $T \subseteq B(\xi, r) \cap [\psi^{-1}(\xi) \setminus (\psi^{\mathrm{GSc}})^{-1}(\xi)]$. But, considering the center-optimal Gale–Shapley algorithm, this implies that $\xi$ rejected $T$, and this contradicts Lemma 19.

Finally, we prove the second inequality in (i). Suppose, on the contrary, that, for every site $x$ in an $\mathcal{L}$-nonnull set $T$ and two fixed centers $\xi, \xi'$, we have $\psi^{\mathrm{GSc}}(x) = \xi$ and $\psi(x) = \xi'$, where $|x - \xi| < |x - \xi'|$. Then either $\xi$ must be unsated in $\psi$, or $\psi^{-1}(\xi)$ must contain a site $y$ further from $\xi$ than a.e. site in $T$, otherwise we would have a contradiction to the first inequality in (ii) just proved. But now for a.e. $x \in T$, we have that the pair $(x, \xi)$ is unstable for $\psi$. □

For the remainder of the article, we take $\Psi = \Psi_\Pi$, where, $\Psi_\bullet$ is any factor allocation, unless otherwise stated. For definiteness, we may take $\Psi_\bullet = \Psi_\bullet^{\mathrm{GS}}$.

## 5. Phases.

PROOF OF THEOREM 4. Using translation-invariance, the result is an immediate corollary of Proposition 20 below. □

Define the *residual appetite* of a center $\xi$ to be $U(\xi) = U_\Psi(\xi) = \alpha - \mathcal{L}\Psi^{-1}(\xi)$.

PROPOSITION 20. *Let $\Pi$ be an ergodic point process of intensity $\lambda \in (0, \infty)$. We have*

$$\mathbf{P}(0 \text{ is unclaimed}) = (1 - \lambda\alpha) \vee 0;$$
$$\mathbf{E}^* U(0) = (\alpha - \lambda^{-1}) \vee 0.$$

PROOF. Lemma 16 with $r = \infty$ states that $\mathbf{P}(0 \text{ is claimed}) = \lambda \mathbf{E}^* \mathcal{L} \Psi_{\Pi^*}^{-1}(0)$. Therefore, from the definition of residual appetite, we have

$$\lambda \mathbf{E}^* U(0) - \mathbf{P}(0 \text{ is unclaimed}) = \lambda\alpha - 1.$$

Recall that no stable allocation may have both unclaimed sites and unsated centers. Also, by ergodicity, the existence of unsated centers and the existence of an $\mathcal{L}$-nonnull set of unclaimed sites are zero–one events. Therefore, at most one of the two terms on the left can be nonzero. The result follows. □

**6. Monotonicity.** We present several results concerning monotonicity in stable allocations. Proposition 21 will be used in the proofs of Theorems 6 and 7. Recall from Section 4 the quantity

$$\gamma_\xi(\psi, r) = \mathcal{L}[\psi^{-1}(\xi) \cap B(\xi, r)],$$

which measures how good an allocation is for a given center.

The following result states that adding extra centers makes sites happier, while making the existing centers less happy.

PROPOSITION 21 (Monotonicity in $\Xi$). *Fix $\alpha$, let $\Xi_1 \subseteq \Xi_2$ be two sets of centers, and let $\psi_1, \psi_2$ be their respective site-optimal Gale–Shapley allocations.*

(i) *For $\mathcal{L}$-a.e. site $x$, we have $|x - \psi_1(x)| \geq |x - \psi_2(x)|$.*
(ii) *For every center $\xi \in \Xi_1$ and all $r \in [0, \infty]$, we have $\gamma_\xi(\psi_1, r) \geq \gamma_\xi(\psi_2, r)$.*

The following result states that increasing the appetite makes sites happier, while making centers more unhappy (in an appropriate sense).

PROPOSITION 22 (Monotonicity in $\alpha$). *Fix $\Xi$, let $\alpha_1 \leq \alpha_2$, and let $\psi_1, \psi_2$ be the respective site-optimal Gale–Shapley allocations to $\Xi$.*

(i) *For $\mathcal{L}$-a.e. site $x$, we have $|x - \psi_1(x)| \geq |x - \psi_2(x)|$.*
(ii) *For every center $\xi \in \Xi$ and all $r \in [0, \infty]$, we have $\alpha_1 - \gamma_\xi(\psi_1, r) \leq \alpha_2 - \gamma_\xi(\psi_2, r)$.*

PROOF OF PROPOSITION 21. (i) We apply Lemma 18 with $\Xi_0 = \Xi_1$, $\Xi = \Xi_2$, $\psi = \psi_1$, $\psi^{\mathrm{GS}} = \psi_2$, $\alpha_0 = \alpha$. The result follows upon recalling that the site-optimal Gale–Shapley algorithm allocates a site to the closest center which does not reject it.

(ii) Suppose the statement fails for some $\xi$ and $r$. Then there exists a positive volume of sites $x \in B(\xi, r)$ with $\psi_2(x) = \xi$ but $\psi_1(x) \notin \{\xi, \Delta\}$, and, furthermore, $\xi$ covets all such sites under $\psi_1$. Pick such an $x$ which is not equidistant from any two centers. By Proposition 21(i), we have $|x - \psi_1(x)| \geq |x - \xi|$, and since $x$ is not equidistant from any two centers, the last inequality is strict. But now $x$ desires $\xi$ under $\psi_1$, so $(x, \xi)$ is an unstable pair for $\psi_1$. □

PROOF OF PROPOSITION 22. (i) We apply Lemma 18 with $\Xi_0 = \Xi$, $\alpha = \alpha_1$, $\alpha_0 = \alpha_2$, $\psi = \psi_1$, $\psi^{\mathrm{GS}} = \psi_2$. The result follows on recalling that the site-optimal Gale–Shapley algorithm allocates a site to the closest center which does not reject it.

(ii) Suppose the statement fails for some $\xi$ and $r$. Then there exists a positive volume of sites $x \in B(\xi, r)$ with $\psi_2(x) = \xi$ but $\psi_1(x) \notin \{\xi, \Delta\}$, and, furthermore, $\xi$ covets all such sites under $\psi_1$. Pick such an $x$ which is not equidistant from any two centers. By Proposition 22(i), we have $|x - \psi_1(x)| \geq |x - \xi|$, and since $x$ is not equidistant from any two centers, the last inequality is strict. But now $x$ desires $\xi$ under $\psi_1$, so $(x, \xi)$ is an unstable pair for $\psi_1$. □

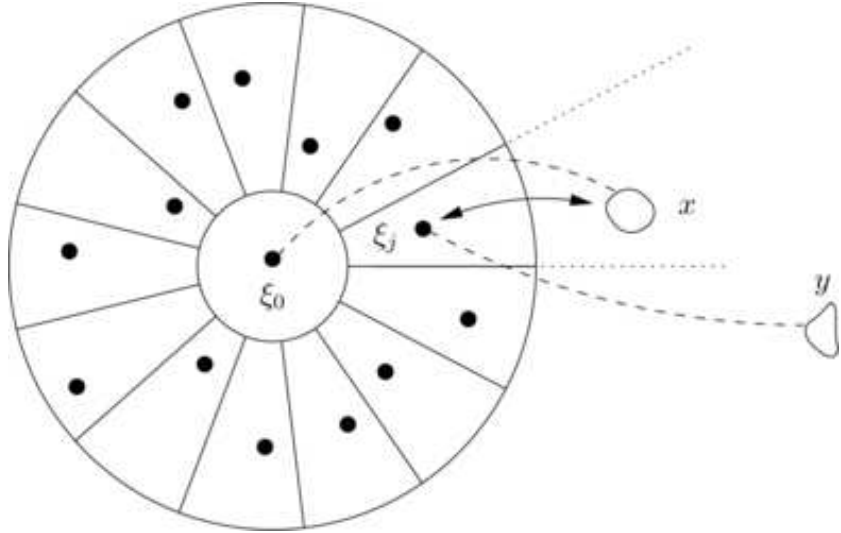

FIG. 3. *An illustration of the proof of Theorem 5(i). If the all centers illustrated have unbounded territory, we obtain a contradiction to stability. Here $(x, \xi_j)$ form an unstable pair.*

**7. Geometry of territories.** In this section we prove Theorem 5. Statement (iii) of that theorem will follow from a somewhat stronger result, Theorem 24, in which we construct a "canonical" version of the stable allocation which has several good properties.

A key tool in proving Theorems 5 and 6 is the following geometric lemma inspired by arguments in [3]. Figure 3 illustrates the main idea. Let $S = \{x \in \mathbb{R}^d : |x| = 1\}$ be the unit sphere. A *cap* is a proper subset of $S$ of the form $H = S \cap B(y, r)$, where $y \in S$. A *cone* is a set of the form

$$V = V_H = \{\alpha h : h \in H \text{ and } \alpha \in (0, \infty)\},$$

where $H$ is a cap.

LEMMA 23 (cones). *For each $d \geq 1$, there exist an integer $k$ and cones $V_1, \ldots, V_k$ whose union is $\mathbb{R}^d \setminus \{0\}$, such that, for any $x, y$ lying in the same cone $V_i$,*

$$\text{if } |x| \geq |y|, \qquad \text{then } |x - y| \leq |x| - \tfrac{1}{2}|y|.$$

(For example, when $d=2$, we may take $k=13$ and $V_i = \{x \neq 0 : \arg x \in [2\pi(i-1)/13, 2\pi i/13)\}$.)

PROOF. The set of all caps of diameter $1/2$ is a cover of the compact set $S$, hence, it has a finite sub-cover $H_1, \ldots, H_k$, say. Put $V_i = V_{H_i}$.

Suppose $x, y \in V_i$ and $|x| \geq |y|$. Let $\beta = |y|/|x|$, so that $|\beta x| = |y|$, and hence, $\beta x, y \in |y| H_i$. Then we have

$$|x - y| \leq |x - \beta x| + |\beta x - y|$$
$$\leq (1-\beta)|x| + \text{diam}(|y|H_i) = |x| - |y| + \tfrac{1}{2}|y|. \quad \square$$

PROOF OF THEOREM 5. (i) Without loss of generality, we may assume that $\Pi$ is ergodic under translations; if not, we may appeal to the ergodic decomposition theorem ([11], Theorem 10.26).

Call a center *bad* if its territory is unbounded. Supposing for a contradiction that bad centers exist, they must form an ergodic translation-invariant point process of positive intensity. Let $V_1, \ldots, V_k$ be cones as in Lemma 23, and note that, since any cone contains arbitrarily large balls, it contains infinitely many bad centers a.s. For $R > 10$, define the sets

$$T_0 = B(0,1);$$
$$T_i = V_i \cap [B(0,R) \setminus B(0,10)], \qquad i = 1, \ldots, k.$$

The above observations imply that if $R$ is chosen sufficiently large, with positive probability, we have

(10) each of $T_0, \ldots, T_k$ contains a bad center.

Suppose the event in (10) occurs, and let $\xi_i$ be a bad center in $T_i$, for each $i$. Since $\xi_0$ is bad, there exists a site $x \in \mathbb{R}^d \setminus B(0,R)$ in its territory. And $x$ lies in some cone, $V_j$, say. Since $\xi_j \in V_j$ also, and $|x| \geq R \geq |\xi_j| \geq 10$, Lemma 23 applies to give

$$|x - \xi_j| \leq |x| - \tfrac{1}{2}|\xi_j|$$
$$\leq |x| - 5$$
$$< |x - \xi_0|.$$

Furthermore, since $\xi_j$ is bad, there exists a site $y$ in its territory with $|y - \xi_j| > |x - \xi_j|$. But then the pair $(x, \xi_j)$ is unstable, a contradiction.

(ii) Without loss of generality, $\Pi$ is ergodic. By translation-invariance, it suffices to prove that $B(0,1)$ intersects only finitely many territories a.s. Note also that, since $\Xi$ is almost surely discrete, if a ball $B$ intersects infinitely many territories, then $\Psi(B)$ must contain centers at arbitrarily large distances from $B$.

Consider the integer lattice $\mathbb{Z}^d$. Call a lattice site $z \in \mathbb{Z}^d$ *bad* if $B(z,1)$ intersects infinitely many territories, and suppose for a contradiction that $\mathbf{P}(0 \text{ is bad}) > 0$. Defining $T_0, \ldots, T_k$ as in the proof of part (i), a similar argument shows that $R$ may be chosen such that, with positive probability,

(11) \quad each of $T_0, \ldots, T_k$ contains a bad lattice site.

On the event in (11), there exists a site $x \in B(0,2)$ with $\Psi(x) \in \mathbb{R}^d \setminus B(R)$, and $\Psi(x) \in V_j$, say. Let $z_j$ be a bad lattice site in $T_j$; then there exists $y \in B(z_j, 1)$ with $|y - \Psi(y)| > |y - \Psi(x)|$. Lemma 23 gives

$$
\begin{aligned}
|y - \Psi(x)| &\leq |z_j - \Psi(x)| + 1 \\
&\leq |\Psi(x)| - \tfrac{1}{2}|z_j| + 1 \\
&\leq |\Psi(x)| - 4 \\
&< |\Psi(x) - x|.
\end{aligned}
$$

So the pair $(y, \Psi(x))$ is unstable. $\square$

The statement in Theorem 5(iii) does not hold for an arbitrary stable factor allocation, since, for instance, we allow $\Psi(x) = \Delta$ on an arbitrary $\mathcal{L}$-null set. Indeed, it does not necessarily hold for $\Psi^{\text{GS}}$, since $\Psi^{\text{GS}}(x) = \Delta$ for all sites which are equidistant from two or more centers. We shall define a "canonical version" $\overline{\Psi}$ of the stable allocation for which Theorem 5(iii) and other desirable properties hold.

Let $\psi$ be a stable allocation, and define a new function $\overline{\psi}$ as follows. For Borel sets $A, B \in \mathbb{R}^d$, we say $A$ is *essentially a subset* of $B$, and write $A \subseteq_{\text{ess}} B$, if $\mathcal{L}(B \setminus A) = 0$. For a site $x$, if there exists $\zeta \in \Xi \cup \{\infty\}$ such that $x$ has some neighborhood which is essentially a subset of $\psi^{-1}(\zeta)$, then we put $\overline{\psi}(x) = \zeta$. If there is no such $\zeta$, we put $\overline{\psi}(x) = \Delta$. Note that $\overline{\psi}$ is well defined because no neighborhood can be essentially a subset of two disjoint sets. Note also that any two allocations which agree $\mathcal{L}$-a.e. yield the same $\overline{\psi}$.

PROOF OF THEOREM 5(iii). Immediate from Theorem 24 below. $\square$

For the remainder of the section, we assume that $\Pi$ is ergodic.

THEOREM 24 (Canonical allocation). *Let $\Psi_\bullet$ be an isometry-equivariant factor allocation to $\Pi$ (e.g., $\Psi_\bullet^{\text{GS}}$) and define $\overline{\Psi}_\bullet$ by the above procedure. Then $\mathbf{P}$-a.s., $\overline{\Psi}_\Pi$ is the unique minimizer of the set $\psi^{-1}(\Delta)$ in the class of stable allocations $\psi$ of $[\Pi]$ having all territories open and the unclaimed set open. Furthermore, $\overline{\Psi}_\bullet$ is an isometry-equivariant factor allocation, and a.s. each territory of $\overline{\Psi}$ has only finitely many connected components, and has $\mathcal{L}$-null boundary.*

In order to prove Theorem 24, we first need to establish some properties of $\Psi$. For a center $\xi$, define
$$\widetilde{R}(\xi) = \sup\{|x - \xi| : x \in \Psi^{-1}(\xi)\}.$$
[Note that $\widetilde{R}(\xi)$ may differ between allocations which agree a.e.] Define the random set of sites
$$Z = \bigcup_\xi \{x \in \mathbb{R}^d : |x - \xi| = \widetilde{R}(\xi)\}$$
$$\cup \bigcup_{\xi \neq \xi'} \{x \in \mathbb{R}^d : |x - \xi| = |x - \xi'|\},$$
where the two unions are over all centers and all pairs of centers respectively. Since $\Xi$ countable, $Z$ is $\mathcal{L}$-null a.s.

LEMMA 25. **P**-a.s., for every site $x \notin Z$, $x$ has a neighborhood $N$ such that $\#(\Psi(N) \setminus \Delta) = 1$.

Before proving Lemma 25, we make some further definitions. Let $B$ be a ball. By Theorem 5(ii), a.s., $B$ intersects only finitely many territories, belonging to centers $\xi_1, \ldots, \xi_k$, say. Define the set
$$Z_B = \bigcup_i \{x \in \mathbb{R}^d : |x - \xi_i| = \widetilde{R}(\xi_i)\}$$
$$\cup \bigcup_{i \neq j} \{x \in \mathbb{R}^d : |x - \xi_i| = |x - \xi_j|\}.$$
Define a *cell* (of $B$) to be any connected component of $B \setminus Z_B$. Since $B \setminus Z_B$ is open, all cells are open sets, and the following lemma implies that, for any given $B$, there are only finitely many cells.

LEMMA 26 (Intersecting hyperplanes). *Let $Y$ be a union of finitely many $(d-1)$-spheres and $(d-1)$-dimensional hyperplanes in $\mathbb{R}^d$. Then $\mathbb{R}^d \setminus Y$ has only finitely many connected components.*

PROOF. By applying an inversion about a site not in $Y$, we may assume that $Y$ is a union of spheres only, and is therefore bounded. Supposing that $\mathbb{R}^d \setminus Y$ has infinitely many components, a repeated bisection argument shows that there must exist a site $z$ all of whose neighborhoods intersect infinitely many components. But considering a neighborhood of $z$ small enough to avoid all those spheres which do not pass through $z$, we see that this is absurd. □

LEMMA 27. **P**-a.s., if $C$ is any cell, then $|\Psi(C) \setminus \Delta| = 1$.

PROOF. We claim first that, for each $x \in B \setminus Z_B$ with $\Psi(x) \neq \Delta$, we have that, $\Psi(x)$ equals the unique closest center to $x$ in the set

$$T(x) = \{\xi_i : |x - \xi_i| < \widetilde{R}(\xi_i)\},$$

or $\Psi(x) = \infty$ if $T(x)$ is empty. To prove this, note first that, by the definitions of $\widetilde{R}(\xi)$ and $Z_B$, we must have $\Psi(x) \in T(x) \cup \{\infty\}$. But if $T(x)$ contains a center $\xi_j$ which is closer to $x$ than $\Psi(x)$ is, then the pair $(x, \xi_j)$ is unstable. This proves the claim.

Now note that the $k^2$ quantities

$$|x - \xi_i| - \widetilde{R}(\xi_i), \quad i = 1, \ldots, k,$$
$$|x - \xi_i| - |x - \xi_j|, \quad i \neq j,$$

are continuous in $x$, and nonzero except on $Z_B$. It follows that their signs are all constant on any given cell. The result now follows from the above claim. □

PROOF OF LEMMA 25. Let $B = B(x, 1)$, and define $Z_B$ and cells as above. Clearly, we have $Z_B \subseteq Z$, so $x$ lies in some cell. Since all cells are open, $x$ has a neighborhood lying in the cell. Now we use Lemma 27. □

PROOF OF THEOREM 24. First note that, since $Z$ is $\mathcal{L}$-null, Lemma 25 implies that a.s. $\overline{\Psi}$ agrees with $\Psi$ $\mathcal{L}$-a.e. It follows immediately that $\overline{\Psi}_\Pi$ is an allocation a.s. By the construction, $\overline{\Psi}_\bullet$ inherits the isometry-equivariance of $\Psi_\bullet^{\mathrm{GS}}$. It also follows immediately from the definition of $\overline{\Psi}$ that the territories and the unclaimed set are open.

To prove that $\overline{\Psi}$ is a stable allocation, suppose, on the contrary, that $(x, \xi)$ is an unstable pair. Suppose first that $\overline{\Psi}(x) = \zeta$ and $\overline{\Psi}(y) = \xi$, where $|x - \xi| < |x - \zeta|$ and $|x - \xi| < |y - \xi|$. Then $x$ and $y$ have neighborhoods intersecting $\Psi^{-1}(\zeta)$ and $\Psi^{-1}(\xi)$ in full measure respectively, so we may find $x', y'$ with $\Psi(x') = \zeta$ and $\Psi(y') = \xi$, and close enough to $x, y$ that we still have $|x' - \xi| < |x' - \zeta|$ and $|x' - \xi| < |y' - \xi|$. Similarly, if $x$ is unclaimed in $\overline{\Psi}$, then we may find a nearby $x'$ which is unclaimed in $\Psi$, and finally, if $\xi$ is unsated in $\overline{\Psi}$, then it is unsated in $\Psi$. In each case we deduce that $(x', \xi)$ is an unstable pair for $\Psi$, a contradiction.

To prove the required minimality statement, we note that a.s., for any $\psi$ as described in the statement of the theorem and any $x$ with $\psi(x) \neq \Delta$, we have that $\psi$ takes the value $\psi(x)$ on some neighborhood of $x$, so $\overline{\Psi}(x) = \psi(x)$.

To prove that the territories of $\overline{\Psi}$ have only finitely many components, note that, in view of Theorem 5(i) and translation-invariance, it is sufficient to rule out the possibility that the intersection of $B(0, 1)$ with some territory has infinitely many components. Let $B = B(0, 1)$ and define $Z_B$ and cells as above. Recall that there are only finitely many cells (by Lemma 26), and that each cell is connected. Lemma 27 combined with the definition $\overline{\Psi}$ implies that $\overline{\Psi}$ is constant on each cell. Hence, it is sufficient to rule out the possibility that the intersection of some territory with $B$ has a component lying entirely in $B \cap Z_B$. But this is impossible since any such component is open, while $Z_B$ is $\mathcal{L}$-null.

Finally, to prove that the territories of $\overline{\Psi}$ have $\mathcal{L}$-null boundaries, note that, by Lemma 25, if $x \notin Z$, then $x$ does not lie in the boundary of any territory. But $Z$ is $\mathcal{L}$-null a.s. $\square$

**8. Infinite demand.** In this section we prove Theorem 6. We shall deduce it from a stronger result which applies to more general point processes.

For a simple point process $\Pi$ and a Borel set $S \subseteq \mathbb{R}^d$ with $\mathcal{L}S \in (0, \infty)$, let $\Pi^S$ be the point process obtained by adding a random extra center in $S$:

$$\Pi^S(A) = \Pi(A) + \delta_U(A),$$

where $U$ is uniformly distributed in $S$, independent of $\Pi$. We call the process $\Pi$ *insertion-tolerant* if the law of $\Pi^S$ is absolutely continuous with respect to the law of $\Pi$ for every $S$. Let $\Pi_S$ be the point process obtained by deleting all centers in $S$:

$$\Pi_S(A) = \Pi(A \setminus S).$$

We call $\Pi$ *deletion-tolerant* if the law of $\Pi_S$ is absolutely continuous with respect to the law of $\Pi$ for every $S$.

THEOREM 28. *Let $\Pi$ be an ergodic point process and let $\lambda \alpha = 1$.*

(i) *If $\Pi$ is insertion-tolerant, then a.s. every center is desired by an infinite volume of sites.*

(ii) *If $\Pi$ is deletion-tolerant, then a.s. every site is coveted by infinitely many centers.*

PROOF OF THEOREM 6. It is easy to check that a Poisson process is both insertion-tolerant and deletion-tolerant, so the result follows immediately from Theorem 28. $\square$

Next we prove Theorem 28. Our starting point is the following weaker result which states that an infinite set of sites would desire some appropriately-positioned center in the unit ball. This result will also be used in the proof of Theorem 7. In the following, we take $|x - \infty| = \infty$ for $x \in \mathbb{R}^d$.

LEMMA 29. *Let $\lambda \alpha = 1$ and let $\Pi$ be an ergodic insertion-tolerant process. We have a.s.*

$$\mathcal{L}\{x \in \mathbb{R}^d : |x - \Psi(x)| \geq |x| - 1\} = \infty.$$

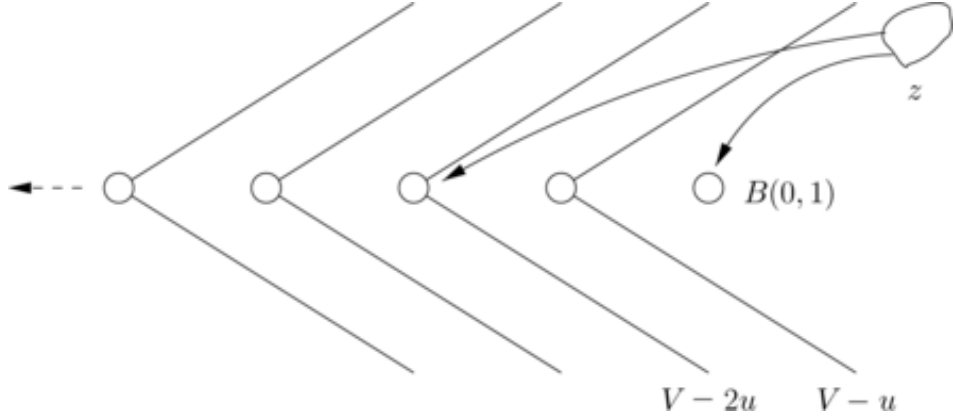

FIG. 4. *An illustration of the proof of Theorem* 28(i). *Some cone must contain an infinite volume of sites which desire all centers in the ball at its apex, and therefore, it must also contain an infinite volume of sites which desire all centers in the ball $B(0,1)$.*

PROOF. Fix $K < \infty$. We shall show that the quantity in question is at least $K$. Let $\Pi'$ be a point process obtained from $\Pi$ by superposing $m$ extra centers, independently uniformly distributed in the ball $B = B(0,1)$, where $m$ is an integer such that $\alpha m \geq K$. By repeatedly applying the insertion-tolerance condition, the law $\mathbf{P}'$ of $\Pi'$ is absolutely continuous with respect to the law $\mathbf{P}$ of $\Pi$. Hence, if $\Psi_\bullet$ is a stable factor allocation as usual, $\Psi_{\Pi'}$ is defined $\mathcal{L}$-a.e. $\mathbf{P}'$-a.s. Also, since $\lambda\alpha = 1$, by Theorem 4, $\mathbf{P}$-a.s., every center is sated in $\Psi_\Pi$. Hence, the same holds $\mathbf{P}'$-a.s. in $\Psi_{\Pi'}$.

The last observation implies that, $\mathbf{P}'$-a.s., we have

$$\mathcal{L}\Psi_{\Pi'}^{-1}(B) \geq \alpha m \geq K. \tag{12}$$

Now since $[\Pi] \subseteq [\Pi']$, Proposition 21(i) gives that, almost surely with respect to the joint law, for $\mathcal{L}$-a.e. $x \in \Psi_{\Pi'}^{-1}(B)$, we have

$$|x - \Psi_\Pi(x)| \geq |x - \Psi_{\Pi'}(x)| \geq |x| - 1.$$

Combining this with (12) yields that, $\mathbf{P}$-a.s.,

$$\mathcal{L}\{x \in \mathbb{R}^d : |x - \Psi_\Pi(x)| \geq |x| - 1\} \geq K. \qquad \square$$

PROOF OF THEOREM 28. (i) The proof is illustrated in Figure 4. Let $V_1, \ldots, V_k$ be the cones of Lemma 23. Since they partition $\mathbb{R}^d \setminus \{0\}$, Lemma 29 implies that, for some fixed $\ell$, we have

$$\mathbf{P}[\mathcal{L}\{x \in V_\ell : |x - \Psi(x)| \geq |x| - 1\} = \infty] > 0. \tag{13}$$

Pick such an $\ell$ and write $V = V_\ell$, and also pick $w \in \mathbb{R}^d$ such that $B(w,1) \subseteq V \setminus B(0,10)$. Note that, for any $n \geq 1$, we also have

(14) $$B(nw, 1) \subseteq V \setminus B(0, 10).$$

For $z \in \mathbb{R}^d$ and $A \subseteq \mathbb{R}^d$, we denote the translated set $A + z = \{a + z : a \in \mathbb{R}^d\}$. Define the set of sites

$$L_n = \{x \in V - nw : |x - \Psi(x)| \geq |x + nw| - 1\},$$

and the event

$$A_n = \{\mathcal{L}(L_n) = \infty\}.$$

The sequence of events $(A_n)_{n \geq 0}$ is stationary and ergodic under $\mathbf{P}$, and (13) implies that $\mathbf{P}(A_0) > 0$. Therefore, we have

(15) $$\mathbf{P}\left[\bigcup_{n=1}^{\infty} A_n\right] = 1.$$

Now suppose that the event $A_m$ occurs, where $m \geq 1$. We shall deduce that every center in $B(0, 1)$ is desired to be an infinite volume of sites. Let

$$z \in L_m \setminus B(-mw, m|w| + 1)$$

(the occurrence of $A_m$ guarantees that an infinite volume of such sites $z$ exists). Let $\xi$ be any center in $B(0, 1)$ (if one exists). We shall apply Lemma 23 with $x = z + mw$ and $y = \xi + mw$. Note that $z \in V - mw$ implies $z + mw \in V$, while (14) implies $\xi + mw \in V$. Furthermore, since $z \notin B(-mw, m|w| + 1)$, we have $|z + mw| \geq m|w| + 1 \geq |mw| + |\xi| \geq |\xi + mw|$. Thus, Lemma 23 applies to yield

$$|z - \xi| \leq |z + mw| - \tfrac{1}{2}|\xi + mw|.$$

Since $z \in L_m$, we have $|z + mw| \leq |z - \Psi(z)| + 1$. By (14), we have $|\xi + mw| \geq 10$. Substituting into the above equation, we obtain

$$|z - \xi| \leq |z - \Psi(z)| + 1 - 5 < |z - \Psi(z)|.$$

Hence, $z$ desires $\xi$.

The above argument together with (15) shows that, almost surely, every center in $B(0, 1)$ is desired by an infinite volume of sites. Since $\mathbb{R}^d$ may be covered by a countable collection of unit balls, by translation invariance, the same claim applies to every center. $\square$

Now we turn to the proof of Theorem 28(ii), which follows similar lines to that of (i). For a center $\xi$, define the *radius* of its territory:

$$R(\xi) = R_\Psi(\xi) = \operatorname*{ess\,sup}_{x \in \Psi^{-1}(\xi)} |\xi - x|.$$

LEMMA 30. *Let $\lambda\alpha = 1$ and let $\Pi$ be an ergodic deletion-tolerant process. We have a.s.*

$$\#\{\xi \in \Xi : R(\xi) \geq |\xi| - 1\} = \infty.$$

PROOF. We mimic the proof of Lemma 29. We show that the quantity in question is at least $K$. Write $B = B(0,1)$. By scaling $\mathbb{R}^d$, without loss of generality, we may assume that $\alpha K \leq \mathcal{L}B$, and the process is still deletion-tolerant. Let $\Pi'$ be the process obtained by deleting all centers in $B$. By Theorem 4 and deletion-tolerance, for $\Pi'$ a.s., infinitely many centers desire some site in $B$. Hence, by Proposition 21(ii), the same is true for $\Pi$. □

PROOF OF THEOREM 28. (ii) We mimic the proof of (i). Take a cone $V$ such that

$$\mathbf{P}(\#\{\xi \in V \cap \Xi : R(\xi) \geq |\xi| - 1\} = \infty) > 0,$$

and take $w$ as before. Almost surely some event

$$A_n = \{\#\{\xi \in (V - nw) \cap \Xi : R(\xi) \geq |\xi + nw| - 1\} = \infty\}$$

occurs, and the previous argument goes through to show that infinitely many centers desire every site in $B(0,1)$. Therefore, the same is true for every site. □

**9. Critical lower bounds.** In this section we deduce Theorem 7 from more general results. The following is Corollary 15 in [10], where it was proved making use of results from [9] and [13].

THEOREM 31. *For $\Pi$ a Poisson process with $\lambda = \alpha = 1$ and for $\Psi_\bullet$ any (not necessarily stable) factor allocation with a.e. site claimed and every center sated, we have $\mathbf{E}X^{d/2} = \infty$ for $d = 1, 2$.*

Interestingly, for $d \geq 3$, there exist "transport rules" (these are a mild generalization of factor allocations) with exponential tails; see [10] for details. We shall prove the following.

THEOREM 32. *Let $\Pi$ be ergodic and insertion-tolerant and let $\lambda\alpha = 1$. For any $d \geq 1$, we have $\mathbf{E}X^d = \infty$.*

PROOF OF THEOREM 7. Immediate from Theorems 31 and 32 above [and Theorem 4(ii)]. □

We prove a somewhat weaker statement for deletion-tolerant processes. Recall the definition of the radius $R(\xi)$ from the last section, and write $R^*(\xi) = R_{\Psi_{\Pi^*}}(\xi)$ for the radius associated with the Palm process.

THEOREM 33. *Let $\Pi$ be ergodic and deletion-tolerant and let $\lambda\alpha = 1$. For any $d \geq 1$, for the Palm process $\Pi^*$, we have $\mathbf{E}^* R^*(0)^d = \infty$.*

PROOF OF THEOREM 32. Since $\Pi$ is insertion-tolerant, Lemma 29 applies. Taking the expectation of the quantity in Lemma 29 and using Fubini's theorem and translation-invariance, we have

$$\mathbf{E}\mathcal{L}\{x : |x - \Psi(x)| \geq |x| - 1\}$$
$$= \int_{\mathbb{R}^d} \mathbf{P}(|x - \Psi(x)| \geq |x| - 1)\mathcal{L}(dx)$$
$$= \int_{\mathbb{R}^d} \mathbf{P}(X \geq |x| - 1)\mathcal{L}(dx).$$

For a constant $c = c(d) \in (0, \infty)$, the last integral may be written

$$\int_{r=0}^{\infty} \mathbf{P}(X \geq r - 1) c r^{d-1} \, dr.$$

Writing $\mu$ for the law of the random variable $X + 1$ and using Fubini's theorem again, the above equals

$$c \int_{z=0}^{\infty} \int_{r=0}^{z} r^{d-1} \, dr \, \mu(dz) = \frac{c}{d} \mathbf{E}[(X+1)^d].$$

Hence, Lemma 29 implies that the above quantity is $\infty$. $\square$

PROOF OF THEOREM 33. The expectation of the quantity in Lemma 30 equals

$$\int_{\mathbb{R}^d} \mathbf{P}^*(R^*(0) \geq |x| - 1)\mathcal{L}(dx);$$

so this quantity is infinite. We now proceed as in the proof of Theorem 32. $\square$

**Open problems.**

(i) *Nonuniqueness.* Does there exist a discrete set of centers in $\mathbb{R}^d$ which admits two stable allocations (differing on a nonnull set)? By Theorems 2 and 3, such a set must have greater than quadratic growth, and may not arise from a translation-invariant point process.

(ii) *Critical tail behavior.* What is the tail behavior of $X$ for the critical Poisson model? In particular, give any quantitative upper bound for $d \geq 2$.

(iii) *Percolation.* In the Poisson model, does the set of unclaimed sites have an infinite component for $\alpha$ sufficiently small? Does it fail to have an infinite component for sufficiently large $\alpha < 1$? Similar questions apply to the set of claimed sites. After this paper was accepted for publication, progress on these problems was achieved in [5].

(iv) *Connected territories.* Is there a translation-equivariant allocation (not stable) of the Poisson process in which every territory is connected, in the critical two-dimensional case ($d=2$ and $\lambda\alpha=1$)?

(v) *Hyperbolic space.* The concept of a stable allocation extends to more general settings, such as hyperbolic space in place of $\mathbb{R}^d$. The proofs of Theorems 1, 3 and 4 adapt to this setting, but the proof of Theorem 5 does not. Are all territories bounded in the setting of hyperbolic space?

**Acknowledgments.** We thank Alan Hammond for suggesting the subcritical and supercritical models, and for valuable conversations. Christopher Hoffman and Yuval Peres thank IMPA in Rio de Janeiro, where some of this work was done. We thank the referee for helpful comments.

C. HOFFMAN  
DEPARTMENT OF MATHEMATICS  
UNIVERSITY OF WASHINGTON  
SEATTLE, WASHINGTON 98195  
USA  
E-MAIL: hoffman@math.washington.edu

A. E. HOLROYD  
DEPARTMENT OF MATHEMATICS  
UNIVERSITY OF BRITISH COLUMBIA  
VANCOUVER, BRITISH COLUMBIA  
CANADA V6T 1Z2  
E-MAIL: holroyd@math.ubc.ca

Y. PERES  
DEPARTMENTS OF STATISTICS  
 AND MATHEMATICS  
UNIVERSITY OF CALIFORNIA, BERKELEY  
BERKELEY, CALIFORNIA 94720  
USA  
E-MAIL: peres@stat.berkeley.edu